\documentclass[english]{ccdconf}
\usepackage{amssymb, amsmath}
\usepackage{array}
\usepackage{multirow}
\usepackage{moreverb}
\usepackage{graphicx,epsfig}
\date{}
\usepackage{caption}
\usepackage{algorithm}
\usepackage{algorithmicx}
\usepackage{algpseudocode}
\usepackage{amsmath}
\usepackage{color}
\usepackage{algpseudocode}
\usepackage{float}
\usepackage{amsthm}
\usepackage{amsfonts}
\usepackage{amssymb}
\usepackage{color}
\usepackage{mathrsfs}
\usepackage{amsmath}
\usepackage{cases}
\begin{document}
\pdfoutput=1

\title{Model-free Value Iteration Algorithm for Continuous-time Stochastic Linear Quadratic Optimal Control Problems}
\author{Guangchen Wang\aref{amss},Heng Zhang\aref{hit}}
\affiliation[amss]{School of Control Science and Engineering,
Shandong University, Jinan, 250061, China.
        \email{wguangchen@sdu.edu.cn}}
\affiliation[hit]{School of Control Science and Engineering,
	Shandong University, Jinan, 250061, China.
        \email{zhangh2828@163.com}}

\maketitle

\begin{abstract}
This paper presents a novel value iteration (VI) algorithm for finding the optimal control for a kind of infinite-horizon  stochastic linear quadratic (SLQ) problem with unknown systems. First, an off-line algorithm is estabilished to obtain the optimal feedback control of our problem. Then, based on the off-line algorithm, the VI-based model-free algorithm and its convergence proof is provided. The main feature of the model-free algorithm is that a stabilizing control is not needed to initiate the algorithm. Finally, we validate our results with a simulation example.
\end{abstract}


\footnotetext{The authors acknowledge the financial support from the NSFC under Grant Nos. 61821004, 11831010 and 61925306, and the NSF of Shandong Province under Grant Nos. ZR2019ZD42 and ZR2020ZD24.}

\section{INTRODUCTION}

Linear quadratic (LQ) control initiated by Kalman \cite{Kalman1960} is an important class of optimal control problems. The SLQ optimal control problem pioneered by Wonham \cite{Wonham1968} is also a fundamental tool in both  application and theory, which has been widely studied by many researchers in the past literature (e.g. Yong and Zhou \cite{YongZhou1999}, Rami et al. \cite{RamiChen2001}, Rami and Zhou \cite{RamiZhou2000}, Zhang, Zhang and Chen \cite{Zhangweihai2008} and Liu, Li and Zhang \cite{LiuLiZhang2014}).
As is known to all, the SLQ optimal control problem in infinite horizon will result in a stochastic algebraic Riccati equation (SARE), which is difficult to solve due to the nonlinear structure of it. With the development of some mathematical methods, many researchers have proposed a lot of methods to get the solution to the corresponding SARE of their problems. See, for instance, Rami and Zhou \cite{RamiZhou2000}, Li, Tai and Tian \cite{LiTaiTian2021} and Liu, Li and Zhang \cite{LiuLiZhang2014}. However, the methods mentioned above are all based on all information of their systems, i.e. all parameters of their systems have to be known beforehand. However, in many pratical problems, the system parameters is hard to obtain. Thus, it is of great importance to solve the SLQ problem without using the information of the system dynamics.\\

In recent years, adaptive dynamic programming (ADP) (Werbos \cite{Werbo1974}) and reinforcement learning (RL) (Sutton and Barto \cite{SuttonBarto1998}) theories have attracted the attention of many researchers. One of the most important applications of RL and ADP is to solve optimal control problems with partially model-free or model-free dynamics. For example, in deterministic system case, one can see Vrabie et al. \cite{Vrabie2009}, Jiang and Jiang \cite{JiangJiang2012}, Bian and Jiang \cite{Bian2016}, Vrabie and Lewis \cite{VrabieNonlinear2009}, Wei, Zhang and Dai \cite{Weizhangdai2009}, Liu et al. \cite{Liu2012} and the references therein. In parallel,  for stochastic optimal control problems,  Wang, Zhang and Luo \cite{WangZhangLuo2016} and Wang, Zhang and Luo \cite{WangZhangLuo2018} got the optimal control for two kinds of linear model-free systems by the method of ADP and Q-learning, respectively. In Ge, Liu and Li \cite{Ge2021}, the authors employ the method of Q-learning in RL to solve the optimal control problem with a class of mean-field discrete-time stochastic systems, without using the information of the system parameters.  Wu and Shen \cite{Wu2018} developed the methodology of policy iteration to obtain the optimal policy of a class of stochastic logical systems.\\

In this work, we will employ the technique of ADP to obtain the optimal policy for a kind of SLQ optimal control problems where the diffusion and drift terms in dynamics depends on both the control and state variables.  We will propose two algorithms to deal with the infinite-horizon SLQ optimal control problem in  continuous time. The first is an off-line algorithm, which borrows the idea of stochastic approximation (e.g. Andrieu, Moulines and  Priouret \cite{Andrieu2005}, Ljung \cite{Ljung1977}, Abounadi, Bertsekas and Borkar \cite{Abouandi2002},  Kushner and Clark \cite{KushnerClark1978} and Kushner and  Yin \cite{KushnerYin2003}). Then, the online VI-based model-free ADP algorithm is presented and the convergence is proved under mild conditions. The proposed online VI-based ADP algorithm has at least two advantages: (1) the proposed algorithm is a  model-free ADP algorithm, i.e. all  system matrices are not required in the process of implementing our algorithm. (2) a stabilizing control policy is not needed to start the algorithm.\\

This work is organized as follows. In section \ref{sec2}, the SLQ optimal control problem and some preliminaries are introduced. In section \ref{sec3}, we propose an off-line algorithm and its convergence is also presented. Based on the first algorithm, a novel VI-based model-free ADP algorithm is presented with rigorous convergence analysis in section \ref{sec4}. In section \ref{sec5}, a numerical example is provided to illustrate the obtained algorithm. Finally, in section \ref{sec6}, some conclusions are given.\\

\noindent{\bf Notations.} We denote by $\mathbb{R}$ the set of real numbers, by $\mathbb{Z^{+}}$ the set of non-negative integers and by $\mathbb{Z^{++}}$ the set of positive integers, respectively. The collection of all $m\times n$ real matrices is denoted by $\mathbb{R}^{m\times n}$. $\mathbb{R}^{m}$ represents the $m$-dimensional Euclidean space and $|\cdot|$ is the Euclidean norm for matrix or vector of proper size. 
The transpose of a vector or matrix $C$ is denoted by $C^{T}$. $\textbf{S}^{m}$, $\textbf{S}^{m}_{+}$ and $\textbf{S}^{m}_{++}$ denote the collection of all symmetric matrices, positive semidefinite matrices and positive definite matrices in $\mathbb{R}^{m\times m}$, respectively. For any matrix $C\in \textbf{S}^{m}_{++}$ (resp. $C\in\textbf{S}^{m}_{+}$), we usually write $C>0$ (resp. $C\geqslant 0$). And for matrices $C\in \textbf{S}^{m}$, $D\in \textbf{S}^{m}$, we write $C>D$ (resp. $C\geqslant D$) if $C-D>0$ (resp. $C-D\geqslant 0$). A function $g$ is belong to $\mathscr{C}(\mathcal{D})$ if $g$ is continuous on $\mathcal{D}$, where $\mathcal{D}\subset\mathbb{R}^{m}$.  
($\Omega$,\,$\mathcal{F}$,\,$\{\mathcal{F}_{t}\}_{t\geqslant 0}$,\,$\mathbb{P}$) represents a filtered probability space that satisfies usual conditions, and $w(\cdot)$ is a one-dimensional standard Brownian motion defined on it.
The Hilbert space is defined by 
\begin{equation*}
	\begin{split}
		L^2_{\mathcal{F}}(\mathbb{R}^{n})=\bigg\{&\psi (\cdot):[0,+\infty)\times \Omega\ \to \mathbb{R}^{n}\bigg|\psi(\cdot) \,\,\text{is}\\ &\mathcal{F}_{t}-\text{adapted,\,\,measureable,\,\,and}\\
		&\,\,  \mathbb{E}\int_{0}^{\infty} |\psi(s,\omega )|^2ds<\infty \bigg\}
	\end{split}	
\end{equation*}
with the corresponding norm  
\begin{equation*}
	\left \|\psi (\cdot) \right \|=\big(\mathbb{E}\int_{0}^{\infty} |\psi(s,\omega )|^2ds\big)^{\frac{1}{2}}.
\end{equation*}

Furthermore,  $\otimes$ is the Kronecker product and for any matrix $C\in\mathbb{R}^{m\times n}$, we define $vec(C)$ as  $vec(C)=[c_1^T,c_2^T,\cdots,c_n^T]$, where $c_i$, $i=1,2,\cdots,n$ is the $i$th column of $C$. 
For any $C\in \textbf{S}^{m}$ and $a_{ij}$ is the $(i,j)$th element of matrix $C$, denote $vecs(C)=[c_{11}, c_{12}, \cdots, c_{1m}, c_{22},c_{23},\cdots, c_{m-1m},c_{mm}]^T$.\\
\
\
\section{PROBLEM FORMULATION AND SOME PRELIMINARIES}\label{sec2}

Consider a stochastic linear system 
\begin{equation}
	\label{eq2}
	\begin{cases}
		\begin{split}
			dx(s)=&[Ax(s)+Bu(s)]ds\\
			&+[Cx(s)+Du(s)]dw(s),
		\end{split}\\
		x(0)=x_0,
	\end{cases}
\end{equation}
where $x_0\in \mathbb{R}^{n}$, $u(\cdot)\in\mathbb{R}^m$ and $A$, $C$, $B$, $D$ are given constant matrices of proper sizes. \\

And the cost functional is
\begin{equation}\label{eq3}
	\begin{split}
		J(u(\cdot))=\mathbb{E}\int_{0 }^{\infty}&[u(s)^TRu(s)+x(s)^TQx(s)]ds,
	\end{split}
\end{equation}
where $R>0$, $Q\geq0$ and  $[A,C|Q]$ is exactly observable.\\

Now we give the following definition, which is crucial for the infinite-horizon SLQ optimal control problem.\\

\noindent{\bf Definition 1} System (\ref{eq2}) is called  mean-square stabilizable if for every $x_0$, there is a matrix $K\in \mathbb{R}^{m\times n} $ such that,  the solution of the following equation
\begin{equation}
	\label{eq4}
	\begin{cases}
		\begin{split}
			dx(s)=&(A+BK)x(s)ds\\
			&+(C+DK)x(s)dw(s),
		\end{split}\\
		x(0)=x_0
	\end{cases}
\end{equation}
satisfies $\lim \limits_{s \to \infty}E[x(s)^Tx(s)]=0$.
In this case,  $u(\cdot)=Kx(\cdot)$ is called a (mean-square) stabilizing control policy.\\

Then we give next assumption to obviate some trivial cases.\\

\noindent{\bf Assumption 1} System (\ref{eq2}) is mean-square stabilizable.\\

 Next we define 
\begin{equation}
	\mathcal{U}_{ad}=\{u(\cdot)\in L^2_{\mathcal{F}}(\mathbb{R}^{m})|u(\cdot) \,\, \text{is\, stabilizing}\}
\end{equation}
as the admissible control set and thus present the SLQ problem as follows:\\

\noindent{\bf Problem (SLQ)} For any given $x_0\in \mathbb{R}^{n}$, our task is to find an optimal control $u^*(\cdot) \in \mathcal{U}_{ad}$ such that
\begin{equation}
	J(u^*(\cdot))=\inf \limits_{u(\cdot) \in \mathcal{U}_{ad}}J(u(\cdot)).
\end{equation}

By Rami and Zhou \cite{RamiZhou2000}, Problem (SLQ) leads to the solution of the following SARE
\begin{equation}\label{SARE}
\begin{cases}
	\begin{split}
A^TP&+PA+Q+C^TPC-(PB+C^TPD)\\
&\times(R+D^TPD)^{-1}(B^TP+D^TPC)=0,
\end{split}\\
R+D^TPD>0.
\end{cases}
\end{equation} 
In addition, $u^*(\cdot)=-(R+D^TP^*D)^{-1}(B^TP^*+D^TP^*C)x(\cdot)$ is an optimal control of Problem (SLQ), where $P^*$ is the maximal solution of SARE (\ref{SARE}).\\
 
Next, we give the following results, which is important for us to prove the off-line algorithm proposed in the next section.\\ 

\noindent{\bf Theorem 1} Assume Assumption 1 holds. Given the stochastic differential Riccati equation (SDRE) as follows

\begin{equation}\label{SDRE}
	\begin{cases}
		\begin{split}
			\dot{P}(t)&+A^TP(t)+Q+P(t)A+C^TP(t)C\\
			&-(P(t)B+C^TP(t)D)(R+D^TP(t)D)^{-1}\\
			&\times (B^TP(t)+D^TP(t)C)=0,
		\end{split}\\
		P(T)=M,\\  
		R+D^TP(t)D>0,
	\end{cases}
\end{equation} then, for any $M\in\textbf{S}^{n}_{+}$, the solution $P(t)$ to the SDRE (\ref{SDRE}) satisfies $\lim_{t\rightarrow -\infty}P(t)=P^*$ and $P(t)$ is monotonically nondecreasing as time $t$ decreases.\\

\noindent{\bf Proof} Since $[A,C|Q]$ is exactly observable, then Theorem 4.1, Theorem 4.2 and Theorem 4.6 of Rami et al. \cite{RamiChen2001} imply that our results hold. This  completes the proof.$\hfill\blacksquare$\\

\noindent{\bf Remark 1.} Since SDRE (\ref{SDRE}) is a backward differential matrix equation, we can reverse the timeline in (\ref{SDRE}) to get a forward differential matrix equation
\begin{equation}\label{nizhuan}	
\begin{cases}
\begin{split}
\dot{P}(t)=&P(t)A+Q+A^TP(t)+C^TP(t)C\\
&-(P(t)B+C^TP(t)D)(R+D^TP(t)D)^{-1}\\
&\times(B^TP(t)+D^TP(t)C),
\end{split}	\\
P(0)=M,\\  
R+D^TP(t)D>0.
\end{cases}
\end{equation}

Obviously, from Theorem 1, we know $\lim_{t\rightarrow \infty}P(t)=P^*$, where $P(t)$ is the solution of (\ref{nizhuan}) with $P(0)=M\in\textbf{S}^{n}_{+}$.\\

\
\
\section{AN OFF-LINE ALGORITHM FOR THE SLQ PROBLEM}\label{sec3}

In this section, we will provide an off-line iterative algorithm and the convergence proof is also presented.\\

\begin{algorithm}[h]
	\caption{}
	\label{}
	\begin{algorithmic}[1]
		\State Choose $P_0>0$, $q\leftarrow 0, k\leftarrow 0$.
		
		\Repeat
		
		\State \begin{equation*}
			\begin{split}
				\widetilde{P}_{k+1}\leftarrow
				&\epsilon_k\big(A^TP_k+C^TP_kC+Q+P_kA\\
				&-(P_kB+C^TP_kD)(R+D^TP_kD)^{-1}\\
				&\times (B^TP_k+D^TP_kC)\big)+P_k
			\end{split}
		\end{equation*}
		
		\If{$\widetilde{P}_{k+1}\in D_q$}  
		\State $P_{k+1}\leftarrow\widetilde{P}_{k+1}$  
		\Else  
		\State $P_{k+1}\leftarrow P_0,q\leftarrow q+1$    
		\EndIf 
		
		\State $k\leftarrow k+1$
		
		\Until{$|\widetilde{P}_{k+1}-P_k|/\epsilon_k<\varepsilon$}
	\end{algorithmic}
\end{algorithm}
 First, we give some definitions as follows. $\{{D_q}\}_{q=0}^\infty$ satisfies

\begin{equation}		
	D_q\subseteq D_{q+1}, \lim\limits_{q\rightarrow \infty}D_q=\textbf{S}^{n}_{+}, q \in \mathbb{Z}^{+}.
\end{equation}
and $D_q$, $q=0,1,2,\cdots,$ is bounded with nonempty interior.  
$\epsilon_k\in\mathbb{R}$ and $\{{\epsilon_k}\}_{k=0}^\infty$  satisfy

\begin{equation}		
	\sum \limits_{k=0}^{\infty}\epsilon_k=\infty, 
	\lim \limits_{k\rightarrow\infty}\epsilon_k=0.
\end{equation}

Our first algorithm is provided in Algorithm 1, whose proof is given in Theorem 2.

\noindent{\bf Theorem 2} Consider $\{{P_k}\}_{k=0}^\infty$ defined in Algorithm 1, we have $\lim_{k\rightarrow \infty}P_k=P^*$.\\

Before proving Theorem 2, we give next lemma to show that $\{{P_k}\}_{k=0}^\infty$ is bounded.\\

\noindent{\bf Lemma 1} There exists a compact set $\kappa$ and $N\geq0$ such that $P^*\in \kappa$ and $\{{P_k}\}_{k=N}^\infty\subset\kappa$.\\

\noindent{\bf Proof} 
First, note that, for any $q\in\mathbb{Z^{+}}$, $vecs(\cdot)$ is an isometric isomorphism from $\textbf{S}^{q}$ to $\mathbb{R}^{q(q+1)/2}$. Thus we can rewrite (\ref{nizhuan}) as 
\begin{equation}\label{ODE}
	\dot{p}=g(p),
\end{equation}
where $p=vecs(P)$, $\mathscr{P}:=\{P\in\textbf{S}^n\,\,\big|\,\,R+D^TPD>0\}$ and $g(\cdot):\mathbb{R}^{n(n+1)/2}\cap vecs(\mathscr{P})\rightarrow\mathbb{R}^{n(n+1)/2}$ is 
\begin{equation}
	\begin{split}
			g(p)=&vecs\big(A^TP+C^TPC+Q+PA\\
		&-(PB+C^TPD)(D^TPD+R)^{-1}\\
		&\times (B^TP+D^TPC)\big).
	\end{split}
\end{equation}
By Remark 1, we can see that if $P(0)\in\textbf{S}^{n}_{+}$, the solution $P(\cdot)$ to the equation (\ref{nizhuan}) converges to $P^*$. Thus, it is easy to see that $vecs(P^*)=p^*$ is a locally asymptotically stable equilibrium of (\ref{ODE}). We denote the region of attraction of $p^*$ by $R_A$ and we know 
\begin{equation*}
	R_A\supset vecs(\textbf{S}^{n}_{+}).
\end{equation*}
Then this lemma can be proved by following the proof of Lemma 3.4 in Bian and Jiang \cite{Bian2016}. This completes the proof.$\hfill\blacksquare$\\

Now, based on Lemma 1, we give the proof of Theorem 2.\\

\noindent{\bf Proof of Theorem 2} Choose $N$ as in Lemma 1, from Algorithm 1, we have
\begin{equation*}
	\begin{split}
		{P}_{k+1}=
		P_k&+\epsilon_k\big(A^TP_k+C^TP_kC+Q+P_kA\\
		&-(P_kB+C^TP_kD)(D^TP_kD+R)^{-1}\\
		&\times (B^TP_k+D^TP_kC)\big)+G_k, \forall k\geq N
	\end{split}
\end{equation*}
and $P_k\in\kappa$, $\forall k\geq N$, where $G_k$ is   
\begin{equation*}
	G_k:=\begin{cases}
		P_0-\widetilde{P}_{k+1}, \text{if}\,\,\, \widetilde{P}_{k+1} \notin \kappa,\\
		0, \quad\quad\quad\quad \text{if}\,\,\, \widetilde{P}_{k+1} \in \kappa.
	\end{cases}
\end{equation*}
Now we define
\begin{equation*}
	P^0(t):=\begin{cases}
		P_k, \,\,\,t\in [t_k,t_{k+1}),\\
		P_0,\,\,\, t\leq 0,\\
	\end{cases}
\end{equation*}
\begin{equation*}
	P^k(t):=P^0(t+t_k),
\end{equation*}
where $t_k=\sum_{i=0}^{k-1}\epsilon_i$ for $k\geq1$, $ t\in(-\infty, \infty)$ and $t_0=0$. By above notations, for all $t\geq 0$ and $k\geq N$, we have
\begin{equation}\label{chazhi}
	\begin{split}
			P^k(t)=&\sum_{i=k}^{n(t+t_k)-1} \epsilon_i\big( A^TP_i+C^TP_iC+Q+P_iA\\
			&-(P_iB+C^TP_iD)(R+D^TP_iD)^{-1}\\
			&\times (B^TP_i+D^TP_iC)  \big)+P_k+\sum_{i=k}^{n(t+t_k)-1}Z_i,\\	
			=&P^k(0)+F^k(t)+E^k(t)+G^k(t),\\	
		\end{split}
\end{equation}
where
\begin{equation*}
	\begin{split}
		F^k(t)=\int_{0}^{t}&\big( A^TP^k(s)+Q+C^TP^k(s)C+P^k(s)A\\
		&-(P^k(s)B+C^TP^k(s)D)\\
		&\times(R+D^TP^k(s)D)^{-1}\\
		&\times (B^TP^k(s)+D^TP^k(s)C)  \big)ds,
	\end{split}
\end{equation*}
\begin{equation*}
	G^k(t)=\sum_{i=k}^{n(t+t_k)-1}Z_i,\,\,\, 
	n(t)=\begin{cases}
			0,\,\,t<0,\\
		i,\,\,0\leq t_i\leq t<t_{i+1},\\
	\end{cases}
\end{equation*}
and 
\begin{equation*}
\begin{split}
	E^k(t)=&\sum_{i=k}^{n(t+t_k)-1} \epsilon_i\big( A^TP_i+C^TP_iC+Q+P_iA\\
	&-(P_iB+C^TP_iD)(R+D^TP_iD)^{-1}\\
	&\times (B^TP_i+D^TP_iC)  \big)-F^k(t).
\end{split}
\end{equation*}

In the above definition, if $t \in[0,\epsilon_k)$,  we assume $\sum_{i=k}^{n(t+t_k)-1}*=0$.\\

Then, following similar procedure in the proof of Theorem 3.3 in Bian and Jiang \cite{Bian2016}, we know $\lim_{k\rightarrow \infty}P_k=P^*$. The proof is completed. $\hfill\blacksquare$\\
\
\
\section{MODEL-FREE VI-BASED ADP ALGORITHM FOR THE SLQ PROBLEM}\label{sec4}

In this section,  on the basis of Algorithm 1, an online VI-based ADP algorithm for Problem (SLQ) and its proof are presented.\\ 

First, we use the Ito's formula to $x(s)P_kx^T(s)$, and from (\ref{eq2}), one gets 
\begin{equation}\label{itoeq}
	\begin{split}
		&d\big(x(s)^TP_kx(s)\big)\\
		=&\big\{x(s)^T\big(A^TP_k+P_kA+C^TP_kC\big)x(s)\\
		&+2u(s)^T\big(B^TP_k+D^TP_kC\big)x(s)\\
		&+u(s)^TD^TP_kDu(s)\big\}ds
		+\big\{\cdots\big\}dw(s).\\
	\end{split}
\end{equation}
Integrating from $t$ to $t+\triangle  t$ and taking expection $\mathbb{E}$ on both sides of (\ref{itoeq}), we have
\begin{equation}\label{E}
	\begin{split}
		&\mathbb{E}\big[x(t+\triangle  t)^TP_kx(t+\triangle  t)-x(t)^TP_kx(t)\big]\\
		=&\mathbb{E}\int_{t}^{t+\triangle  t}x(s)^TM_kx(s)
		+2u(s)^TN_kx(s)ds\\
		&+\mathbb{E}\int_{t}^{t+\triangle  t}u(s)^TH_ku(s)ds,\\
		\end{split}
	\end{equation}
where $M_k=A^TP_k+P_kA+C^TP_kC\in\mathbb{R}^{n\times n}$, $N_k=B^TP_k+D^TP_kC\in\mathbb{R}^{m\times n}$ and $H_k=D^TP_kD\in\mathbb{R}^{m\times m}$. Then, for any $\xi\in\mathbb{R}^q$, we define
\begin{equation}
	\overline{\xi}=[\xi_1^2,2\xi_1\xi_2, \cdots, 2\xi_1\xi_q,\xi_2^2,2\xi_2\xi_3, \cdots,2\xi_{q-1}\xi_q, \xi_q^2]^T,
\end{equation}
where $\xi_i$ is the $i$th element of $\xi$ and $q\in\mathbb{Z^{++}}$. Moreover, we define matrices $I _{xx} $, $\delta _{xx} $, $\delta _{xu} $, $\delta _{uu}$ as follows
\begin{equation*}
	\begin{split}
		I _{xx}=\mathbb{E}\big[\bar x(t_1)-\bar x(t_0),\cdots,
		\bar x(t_l)-\bar x(t_{l-1})\big]^T,\\
	\end{split}
\end{equation*}
\begin{equation*}
	\begin{split}
		\delta _{xx}=\mathbb{E}\bigg[\int_{t_0}^{t_1} \overline{x(s)}ds,\cdots,
		\int_{t_{l-1}}^{t_l}\overline{x(s)}ds\bigg]^T,\\
	\end{split}
\end{equation*}
\begin{equation*}
	\begin{split}
		\delta _{xu}=\mathbb{E}\bigg[\int_{t_0}^{t_1}x(s)\otimes u(s)ds,\cdots,
		\int_{t_{l-1}}^{t_l}x(s)\otimes u(s)ds\bigg]^T,\\
	\end{split}
\end{equation*}
\begin{equation*}
	\begin{split}
		\delta _{uu}=\mathbb{E}\bigg[\int_{t_0}^{t_1}\overline{u(s)}ds,\cdots,
		\int_{t_{l-1}}^{t_l}\overline{u(s)}ds\bigg]^T.\\
	\end{split}
\end{equation*}

By (\ref{E}), we have
\begin{equation}\label{eq12}
	\big[\delta_{xx},2\delta_{xu},\delta_{uu}\big]\begin{bmatrix}
			vecs(M_k)\\
			vec(N_k)\\
			vecs(H_k)
	\end{bmatrix}=I_{xx}vecs(P_k)
\end{equation}
Notice that if $u(\cdot)$ is chosen to satisfy 
\begin{equation}\label{noise}
	rank([\delta_{xx},\delta_{xu},\delta_{uu}])=mn+\frac{n(n+1)}{2}+\frac{m(m+1)}{2},
\end{equation}
 then $\big[\delta_{xx},2\delta_{xu},\delta_{uu}\big]\in\mathbb{R}^{l\times \big(mn+\frac{m(m+1)}{2}+\frac{n(n+1)}{2}\big)}$ has full column rank, and, from (\ref{eq12}), one can easily gets
\begin{equation}\label{solve}
	\begin{bmatrix}
		vecs(M_k)\\
		vec(N_k)\\
		vecs(H_k)\\
	\end{bmatrix}=\Theta \,\,vecs(P_k),
\end{equation}
where
\begin{equation*}
	\begin{split}
		\Theta=&\bigg(\big[\delta_{xx},2\delta_{xu},\delta_{uu}\big]^T\big[\delta_{xx},2\delta_{xu},\delta_{uu}\big]\bigg)^{-1}\\
		&\times\big[\delta_{xx},2\delta_{xu},\delta_{uu}\big]^TI_{xx}.
	\end{split}
\end{equation*} \\

Based on (\ref{solve}), the model-free algorithm is summarized in Algorithm 2.\\

\begin{algorithm}[h]
	\caption{}
	\label{}
	\begin{algorithmic}[1]
		\State Choose $P_0>0$, $q\leftarrow 0, k\leftarrow 0$.
		       Employ $u(\cdot)$ that satisfies (\ref
		       {noise}) as the input to $(\ref{eq2})$ and collect the online data to compute $\Theta$.
		
		\Repeat
		
		\State Solve $(M_k,N_k,H_k)$ from (\ref{solve}) 
			\begin{equation*}
			\begin{split}
				\widetilde{P}_{k+1}\leftarrow
				P_k+\epsilon_k&\big(M_k+Q-N_k^T
				(R+H_k)^{-1}N_k\big)
			\end{split}
		\end{equation*}
		
		\If{$\widetilde{P}_{k+1}\in D_q$}  
		\State $P_{k+1}\leftarrow\widetilde{P}_{k+1}$  
		\Else  
		\State $P_{k+1}\leftarrow P_0,q\leftarrow q+1$    
		\EndIf 
		
		\State $k\leftarrow k+1$
		
		\Until{$|\widetilde{P}_{k+1}-P_k|/\epsilon_k<\varepsilon$}
	\end{algorithmic}
\end{algorithm}

\noindent{\bf Theorem 3} When (\ref{noise}) is satisfied, we have $\lim_{k\rightarrow\infty}P_k=P^*$, where $\{P_k\}_{k=0}^\infty$ are given by Algorithm 2.\\

\noindent{\bf Proof} Since (\ref{noise}) is satisfied,  $\big[\delta_{xx},2\delta_{xu},\delta_{uu}\big]$ has full column rank, which implies that (\ref{solve}) has unique solution. Moreover, by (\ref{E})-(\ref{solve}), $\big(vecs(A^TP_k+P_kA+C^TP_kC), vec(B^TP_k+D^TP_kC),vecs(D^TP_kD)\big)$ is a solution to (\ref{solve}), thus the solution $(M_k,N_k,H_k)$ obtained from (\ref{solve}) is equivalent to $(A^TP_k+P_kA+C^TP_kC, B^TP_k+D^TP_kC,D^TP_kD)$. 
Otherwise, there exists a different solution $(M_k^{'},N_k^{'},H_k^{'})$ of (\ref{E}), and is also the solution of (\ref{solve}). 
Therefore, $\widetilde{P}_{k+1}$ and $P_{k+1}$ obtained from Algorithm 2 is equivalent to the ones in Algorithm 1. Then by Theorem 2, one gets $\lim_{k\rightarrow\infty}P_k=P^*$. This completes the proof. $\hfill\blacksquare$\\
\
\
\section{Numerical example}\label{sec5}

In this section, we will present a numerical example of Algorithm 2. Consider the system coefficients 
 \begin{equation*}
	A=
	\begin{bmatrix}
		
		0  & -0.6\\
		0.6   & -0.3
	\end{bmatrix},
	B=
	\begin{bmatrix}
		
		0.05\\
		0.01
	\end{bmatrix},
\end{equation*} 
\begin{equation*}
	C=
	\begin{bmatrix}
		-0.02 &   0.03\\
		-0.05    &0.02
	\end{bmatrix},
	D=
	\begin{bmatrix}
		0.001\\
		0.03    
	\end{bmatrix},
\end{equation*} 
and the initial state is $x_0=[0.5,-0.1]^T$. We choose $R=1$ and $Q=diag\{0.05,0.1\}$ in the cost functional. \\

When the stoping criterion $|\widetilde{P}_{k+1}-P_k|/\epsilon_k<\varepsilon=10^{-5}$ is satisfied, by Algorithm 2, the approximate solution $\bar{P}^*$ and $\bar{K}^*$ are obtained as
\begin{equation*}
	\bar{P}^*=
	\begin{bmatrix}
		0.2722091 & -0.0427624\\
		-0.0427624    &0.2505643
	\end{bmatrix},
\end{equation*}
\begin{equation*}
	\bar{K}^*=
	\begin{bmatrix}
		-0.0134984 & 0.0298522\\
	\end{bmatrix}.
\end{equation*}

Then, to check the error of our algorithm, we denote
\begin{equation*}
	\begin{split}
		\mathcal{R}_1(P)=&A^TP+Q+C^TPC+PA\\
		&-(PB+C^TPD)(R+D^TPD)^{-1}\\
		&\times(B^TP+D^TPC),
	\end{split}
\end{equation*}
\begin{equation*}
	\begin{split}
		\mathcal{R}_2(P, K)=&(A+BK)^TP+P(A+BK)\\
		&+(C+DK)^TP(C+DK)\\
		&+K^TRK+Q,
	\end{split}
\end{equation*}
and we have
\begin{equation*}
	\mathcal{R}(\bar{P}^*)=
	\begin{bmatrix}
		0.0008297 & -0.0004970\\
		-0.0004970    &0.0012700
	\end{bmatrix},
\end{equation*}
\begin{equation*}
	\mathcal{R}(\bar{P}^*,\bar{K}^*)=
	\begin{bmatrix}
		-0.0008292&   -0.0005174\\
		-0.0005174&    0.0021906
	\end{bmatrix}.
\end{equation*}

Obviously, the errors are almost of size $10^{-4}$. 

\section{CONCLUSION}\label{sec6}

This paper is concerned with an continuous-time SLQ optimal control problem in infinite horizon. We first propose an off-line algorithm to get the maximal solution of the corresponding SARE.  Then, by the off-line algorithm, the VI-based model-free ADP algorithm is presented. Without knowing the knowledge of all system coefficient matrices, this algorithm updates by using the  input and state information collected online. Moreover, our algorithm does not need an stabilizing policy to initiate the algorithm, so the the difficulty in searching for an initial stabilizing control has been overcome. Finally, a numerical example is presented to illustrate the obtained algorithm. It is an interesting topic for us to extend this algorithm
to the case that the control weighting matrix in the cost functional to be indefinite.\\

\end{document}